\documentclass[12pt]{amsart}
\usepackage{amsfonts,amsmath,enumerate}
\usepackage{amssymb,amscd,amsthm}
\usepackage[all]{xy}
\usepackage[pdftex]{graphicx}
\theoremstyle{plain}
\newtheorem{thm}[subsection]{Theorem}

\newtheorem{prop}[subsection]{Proposition}
\newtheorem{cor}[subsection]{Corollary}

\theoremstyle{definition}
\newtheorem{rk}[subsection]{Remark}

\newtheorem{ex}[subsection]{Example}

\numberwithin{equation}{section}
\newcommand{\C}{\mathbb{C}}

\newcommand{\PP}{\mathbb{P}}

\numberwithin{equation}{section}

\begin{document}

\title [On Hodge theory of singular plane curves]
{On Hodge theory of singular plane curves }

\author[Nancy Abdallah]{Nancy Abdallah}
\address{Univ. Nice Sophia Antipolis, CNRS,  LJAD, UMR 7351, 06100 Nice, France.}
\email{Nancy.ABDALLAH@unice.fr}

\subjclass[2000]{Primary  32S35, 32S22; Secondary  14H50}

\keywords{plane curves, Hodge  and Pole order filtrations }

\begin{abstract} The dimensions of the graded quotients of the cohomology of a plane curve complement $U=\PP^2 \setminus C$ with respect to the Hodge filtration are described in terms of simple geometrical invariants. The case of curves with ordinary singularities is discussed in detail. We also give a precise numerical estimate for the difference between the Hodge filtration and the pole order filtration on $H^2(U,\C)$.

\end{abstract}

\maketitle

\section{Introduction}
The Hodge theory of the complement of projective hypersurfaces have received a lot of attention, see for instance Griffiths \cite{G} in the smooth case, Dimca-Saito \cite{DS} and Sernesi \cite{Se} in the singular case. In this paper we consider the case of plane curves and continue the study initiated by Dimca-Sticlaru \cite{DSt} in the nodal case and the author \cite{A} in the case of plane curves with ordinary singularities of multiplicity up to $3$.

In the second section we compute the Hodge-Deligne polynomial of a plane curve $C$, the irreducible case in Proposition \ref{propirred} and the reducible case in Proposition \ref{propred}.
Using this we determine the Hodge-Deligne polynomial of $U=\PP^2 \setminus C$ and then
we deduce in Theorem \ref{mhnU} the dimensions of the graded quotients of $H^2(U)$  with respect to the Hodge filtration.

In section three we consider the case of arrangements of curves having ordinary singularities and intersecting transversely at smooth points and obtain a formula in Theorem   \ref{mhnU2} generalizing the formulas obtained in \cite{DSt}  and in \cite{A} (for this type of curves).
In fact, the results in \cite{A} show that this formula holds in the more general case of plane curves with ordinary singularities of multiplicity up to $3$ (without assuming transverse intersection).

In the forth section we show that the case of plane curves with ordinary singularities of multiplicity up to $4$ (without assuming transverse intersection) is definitely more complicated and the formula in Theorem   \ref{mhnU2} has to be replaced by the formula in Theorem   \ref{mhnUformult4} containing a correction term coming from triple points on one component through which another component of $C$ passes.

In the final section we give some applications, we hope of general interest, expressing  the difference between the Hodge filtration and the pole order filtration on $H^2(U,\C)$ in terms of numerical invariants easy to compute in given situations, see Theorem \ref{F=P} and its corollaries. One example involving a free divisor concludes this note.

\section{Hodge Theory of plane curve complements}

For the general theory of mixed Hodge structures we refer to \cite{De} and \cite{Voisin}.
Recall the definition of the Hodge-Deligne polynomial of a quasi-projective complex variety $X$
$$P(X)(u,v)=\sum_{p,q} E^{p,q}(X)u^pv^q$$
where $E^{p,q}(X)=\sum_s (-1)^sh^{p,q}(H_c^s(X))$,  with  $h^{p,q}(H^s_c(X)) =\dim Gr_F^p Gr_{p+q}^WH^s_c(X,\mathbb{C})$, the mixed Hodge numbers of $H_c^s(X)$.

This polynomial  is additive with respect to constructible partitions, i.e. $P(X)=P(X\setminus Y)+P(Y)$ for a closed subvariety $Y$ of $X$. In this section we determine $P(C)$ for a (reduced) plane curve $C$.

Suppose first that the curve $C$ is irreducible, of degree $N$. Denote by $a_k$, $k=1, ... ,p$ the singular points of $C$, and let $r(C,a_k)$ be the number of irreducible branches of the germ $(C,a_k)$. Let  $\nu: \tilde{C} \rightarrow C$ be the normalization mapping.
Using the normalization map $\nu$ and the additivity of the Hodge-Deligne polynomial, it follows that,
$$
P(C)
=P(C\backslash(C)_{sing})+P((C)_{sing})=P(\tilde{C}\backslash (\cup_k \nu^{-1}(a_k))+p=$$
$$=P(\tilde{C})-\sum_k P(  \nu^{-1}(a_k))+p
= uv-gu-gv+1-\sum_k (r(C,a_k)-1).
$$
Indeed, it is known that for the smooth curve $\tilde C$, the genus $g=g(\tilde C)$ is
exactly the Hodge number $h^{1,0}(\tilde C)=h^{0,1}(\tilde C)$.
Moreover, it is known that one has the formula
\begin{equation}\label{genus-delta}
g=\frac{(N-1)(N-2)}{2} -\sum_k\delta(C,a_k),
\end{equation}
relating the genus, the degree and the local singularities of $C$, and the $\delta$-invariants can be computed using the formula
\begin{equation}\label{Milnor}
2\delta(C,a_k)=\mu(C,a_k)+r(C,a_k)-1,
\end{equation}
where $\mu(C,a_k)$ is the Milnor number of the singularity $(C,a_k)$. For both formulas above, see Milnor, p. 85. This proves the following result.

\begin{prop}\label{propirred} 
With the above notation and assumptions, we have the following for an irreducible plane curve $C \subset \PP^2$.

\begin{enumerate}[(i)]

\item The Hodge-Deligne polynomial of $C$ is given by
$$P(C)(u,v)=uv-gu-gv+1-\sum_k (r(C,a_k)-1),
$$
with $g$ given by the formula \eqref{genus-delta}.
\item $H^0(C)=\C$ is pure of type $(0,0)$.

\item  $H^2(C)=\C$ is pure of type $(1,1)$.

\item  The mixed Hodge numbers of the MHS on $H^1(C)$ are given by
$$h^{0,0}(H^1(C))=\sum_k (r(C,a_k)-1), \  \ h^{1,0}(H^1(C))=  h^{0,1}(H^1(C))= g.$$
In particular, one has the following formulas for the first Betti number of $C$.
$$b_1(C)= \sum_k (r(C,a_k)-1) +2g= (N-1)(N-2)-\sum_k \mu(C,a_k).$$
\end{enumerate}
\end{prop}

Now we consider the case of a curve $C$ having several irreducible components. More precisely,
let $C=\bigcup_{j=1,r}C_j$ be the decomposition of $C$ as a union of irreducible components $C_j$, let $\nu_j: \tilde{C}_j \rightarrow C_j$ be the normalization mappings and set $g_j=g(\tilde{C}_j)$.
Suppose that the curve $C_j$ has degree $N_j$, denote by $a_k^j$ for $k=1,...,p_j$ be the singular points of $C_j$ and let $r(C_j,a_k^j)$ be the number of branches of the germ $(C_j,a_k^j)$. Then the formulas \eqref{genus-delta} and \eqref{Milnor} can be applied to each irreducible curve $C_j$, as well as Proposition \ref{propirred}.

Let $A$ be the union of the singular sets of the curves $C_j$. Let $B$ be the set of points in $C$ sitting on at least two distinct components $C_i$ and $C_j$. For $b \in B$, let $n(b)$ be the number of irreducible components $C_j$ passing through $b$. By definition, $n(b) \geq 2$. Moreover, note that the sets $A$ and $B$ are not disjoint in general, and their union is precisely the singular set of $C$.

Using the additivity of Hodge-Deligne polynomials we get
$$
P(C)=P(C_1\cup\cdots\cup C_r)
=\sum_{j=1}^r P(C_j)+(-1)^{l-1}\sum_{0\leq i_1<\cdots<i_l\leq r}P(C_{i_1}\cap\cdots \cap C_{i_l}).$$
The first sum is easy to determine using Proposition \ref{propirred}.
$$\sum_{j=1}^r P(C_j)
=ruv - \left(\sum_{j=1}^r g_j\right)u - \left(\sum_{j=1}^r g_j\right)v+r - \sum_{j,k}( (r(C_j,a_k^j)-1).$$
Consider now the alternated sum, where $l \geq 2$. The only points of $C$ that give a contribution to this sum are the points in $B$. Now, for a point $b \in B$, its contribution to the alternated sum is clearly given by 
$$c(b)=-{n(b) \choose 2}+ {n(b) \choose 3} -...+(-1)^{n(b)-1}{n(b) \choose n(b)} =-n(b)+1.     $$

\begin{prop}\label{propred} 
With the above notation and assumptions, we have the following for a reducible plane curve $C=\bigcup_{j=1,r}C_j$.

\begin{enumerate}[(i)]

\item The Hodge-Deligne polynomial of $C$ is given by
$$P(C)(u,v)=ruv - \left(\sum_{j=1}^r g_j\right)u - \left(\sum_{j=1}^r g_j\right)v+r - \sum_{j,k}( (r(C_j,a_k^j)-1)- \sum_{b\in B}(n(b)-1).
$$
with $g_j$ given by the formula \eqref{genus-delta}.
\item $H^0(C)=\C$ is pure of type $(0,0)$.

\item  $H^2(C)=\C^r$ is pure of type $(1,1)$.

\item  The mixed Hodge numbers of the MHS on $H^1(C)$ are given by
$$h^{0,0}(H^1(C))=   \sum_{j,k}( (r(C_j,a_k^j)-1)+ \sum_{b\in B}(n(b)-1)   -r+1, $$      
$$ h^{1,0}(H^1(C))=  h^{0,1}(H^1(C))= \sum_jg_j.$$
In particular, one has the following formula for the first Betti number of $C$.
$$b_1(C)= \sum_{j,k}( (r(C_j,a_k^j)-1)+ \sum_{b\in B}(n(b)-1)   -r+1   
+2 \sum_jg_j.$$ 
\end{enumerate}
\end{prop}

Note that a point in the intersection $A \cap B$ will give a contribution to the last two sums in the above formula for $P(C)$.

\begin{ex}\label{exnodes} Suppose $C$ is a nodal curve. Then for each singularity $a_k^j\in A$ one has $a_k^j \notin B$ (otherwise we get worse singularities than nodes) and $r(a_k^j)=2$.
Moreover, each point $b \in B$ satisfies $n(b)=2$. It follows that in this case we get
$$P(C)(u,v)=ruv - \left(\sum_{j=1}^r g_j\right)u - \left(\sum_{j=1}^r g_j\right)v+r -n_2,
$$
with $n_2$ the number of nodes of $C$. More precisely, in this case we have
$n_2=n_2'+n_2''$, where $n_2'$ (resp. $n_2''$) is the number of nodes of $C$ in $A$ (resp. in $B$) and one clearly has
$$n_2'=S_1:=\sum_{j,k}( (r(C_j,a_k^j)-1), \  \ n_2''=S_2:=\sum_{b\in B}(n(b)-1).$$

\end{ex}

\begin{ex}\label{exnodes+triple} Suppose $C$ has only nodes and ordinary triple points as singularities. Then let $n_3$ be the number of triple points and note that we can write as above
$n_3=n_3'+n_3''$, where $n_3'$ (resp. $n_3''$) is the number of triple points of $C$ in $A_0=A\setminus B$ (resp. in $B$). For a point $a \in A_0$, the contribution to the sum $S_1$ is $2$, while the contribution to the sum $S_2$ is $0$.

A point $b \in B$ can be of two types. The first type, corresponding to the partition $3=1+1+1$, is when $b$ is the intersection of three components $C_j$, all smooth at $b$. The contribution of such a point $b$ is $0$ to the sum $S_1$ and $2$ to the sum $S_2$.

The second type, corresponding to the partition $3=2+1$, is when $b$ is the intersection of two components, say $C_i$ and $C_j$, such that $C_i$ has a node at $b$, and $C_j$ is smooth at $b$. The contribution of such a point $b$ is $1$ to the sum $S_1$ and $1$ to the sum $S_2$.

It follows that the contribution of any triple point to the sum $S_1+S_2$ is equal to $2$. Since the double points in $C$ can be treated exactly as in 
 Example \ref{exnodes}, this yields the following.
$$P(C)(u,v)=ruv - \left(\sum_{j=1}^r g_j\right)u - \left(\sum_{j=1}^r g_j\right)v+r -n_2 -2n_3.
$$
When there are only triple points in $B$ of the first type, then we obviously have the following additional relations
$$S_1=n_2'+2n_3', \  \  S_2=n_2''+2n_3''.$$

\end{ex}

\begin{ex}\label{exnodes+triple+4fold} Suppose $C$ has only  ordinary points of multiplicity 2, 3 and 4 as singularities.  Then let $n_4$ be the number of  points of multiplicity 4 and note that we can write as above
$n_4=n_4'+n_4''$, where $n_4'$ (resp. $n_4''$) is the number of  points of multiplicity 4 of $C$ in $A_0=A\setminus B$ (resp. in $B$).
For a point $a \in A_0$ of multiplicity $4$, the contribution to the sum $S_1$ is $3$, while the contribution to the sum $S_2$ is $0$.

A point $b \in B$ can be of 4 types. The first type, corresponding to the partition $4=1+1+1+1$, is when $b$ is the intersection of 4 components $C_j$, all smooth at $b$. The contribution of such a point $b$ is $0$ to the sum $S_1$ and $3$ to the sum $S_2$.

The second type, corresponding to the partition $4=2+1+1$, is when $b$ is the intersection of 3 components, say $C_i$, $C_j$ and $C_k$, such that $C_i$ has a node at $b$, and $C_j$ and $C_k$ are smooth at $b$. The contribution of such a point $b$ is $1$ to the sum $S_1$ and $2$ to the sum $S_2$.

The third type, corresponding to the partition $4=2+2$, is when $b$ is the intersection of 2 components, say $C_i$ and $C_k$, such that $C_i$ and $C_k$ have a node at $b$. The contribution of such a point $b$ is $2$ to the sum $S_1$ and $1$ to the sum $S_2$.

The fourth type, corresponding to the partition $4=3+1$, is when $b$ is the intersection of 2 components, say $C_i$ and $C_k$, such that $C_i$ has a triple point at $b$, and  $C_k$ is smooth at $b$. The contribution of such a point $b$ is $2$ to the sum $S_1$ and $1$ to the sum $S_2$.

It follows that the contribution of any  point of multiplicity 4 to the sum $S_1+S_2$ is equal to $3$. Since the double and triple points in $C$ can be treated exactly as in 
 Example \ref{exnodes+triple}, this yields the following.
$$P(C)(u,v)=ruv - \left(\sum_{j=1}^r g_j\right)u - \left(\sum_{j=1}^r g_j\right)v+r -n_2 -2n_3-3n_4.
$$
When there are only  points of multiplicity 4 in $B$ of the first type, then we obviously have the following additional relations
$$S_1=n_2'+2n_3'+3n_4'' , \  \  S_2=n_2''+2n_3''+3n_4''.$$

\end{ex}

Let's look now at the cohomology of the smooth surface $U=\PP^2 \setminus C$. 
By the additivity we get $P(U)=P(\mathbb{P}^2)-P(C)$ where $P(\mathbb{P}^2)=u^2v^2+uv+1$.
This yields the following consequence.

\begin{cor}\label{PU}
$$P(U)(u,v)=u^2v^2-(r-1)uv + \left(\sum_{j=1}^r g_j\right)u + \left(\sum_{j=1}^r g_j\right)v-(r-1) + $$
$$+\sum_{j,k}( (r(C_j,a_k^j)-1)+ \sum_{b\in B}(n(b)-1).$$

\end{cor}

The contribution of $H^4_c(U,\C)$ to $P(U)$ is the term $u^2v^2$, and that of $H^3_c(U,\C)$ is the term $-(r-1)uv$. Moreover, the dimension $\dim Gr^1_FH^2(U,\C)$ is the number of independent classes of type (1,2), which correspond to classes of type $(1,0)$ in $H^2_c(U)$, and hence to the terms in $u$ in $P(U)$. For both statements see the proof of Theorem 2.1 in \cite{A}. This proves the following result.

\begin{thm} \label{mhnU}
$$\dim Gr^1_FH^2(U,\C)=\sum_{j=1}^r g_j  $$
and
$$\dim Gr^2_FH^2(U,\C)=\sum_{j=1}^r g_j +   \sum_{j,k}( (r(C_j,a_k^j)-1)+ \sum_{b\in B}(n(b)-1)   -r+1. $$
In particular, all the components $C_j$ of the curve $C$ are rational if and only if $H^2(U)$ is pure of type $(2,2)$.
\end{thm}

\begin{ex} \label{exnodes+triple+4fold2} Suppose $C$ has only  ordinary points of multiplicity 2, 3 and 4 as singularities.  Then let $n_k$ be the number of  points of multiplicity $k$, for $k=2,3,4$
Then using Example \ref{exnodes+triple+4fold}, we get the formula
$$\dim Gr^2_FH^2(U,\C)=\sum_{j=1}^r g_j -r+1+  n_2+2n_3+3n_4. $$

\end{ex}

\section{Arrangements of transversely intersecting curves}

Recall that  $C=\bigcup_{j=1,r}C_j$ is the decomposition of $C$ as a union of irreducible components $C_j$,
and the curve $C_j$ has degree $N_j$. 
In this section we assume that any curve $C_j$ has only ordinary multiple points as singularities and let $n_k(C_j)$ denote the number of ordinary points on $C_j$ of multiplicity $k$.
We also assume that the intersection of any two distinct components $C_i$ and $C_j$ is transverse, i.e. the points in $C_i \cap C_j$ are nodes of the curve $C_i \cup C_j$.
This implies in particular that $A \cap B=\emptyset$. The formulas \eqref{genus-delta} and 
\eqref{Milnor} yield the equality.
\begin{equation}\label{genus-delta2}
g_j=\frac{(N_j-1)(N_j-2)}{2} -\frac{1}{2} \sum_k\left( \mu(C_j,a_k^j)+r(C,a_k^j)-1\right),
\end{equation}
Using this, Theorem \ref{mhnU} gives the formula
$$
\dim Gr^2_FH^2(U,\C)=\sum_{j=1}^r\frac{(N_j-1)(N_j-2)}{2} -\frac{1}{2} \sum_{j,k}\left( \mu(C_j,a_k^j)-r(C,a_k^j)+1\right) +$$
$$+\sum_{b\in B}(n(b)-1)   -r+1.
$$
If $a_k^j$ is an ordinary $m$-multiple point on the curve $C_j$, one has $\mu(C_j,a_k^j)=(m-1)^2$ and hence 
$$\mu(C_j,a_k^j)-r(C,a_k^j)+1=(m-1)(m-2).$$
If we denote by $n_m'$ (resp. $n_m''$) the number of $m$-multiple points of $C$ coming from just one component $C_j$ (resp. from the intersection of several components $C_j$), we see that
we have
$$\sum_{j,k}\left( \mu(C_j,a_k^j)-r(C,a_k^j)+1\right)= \sum_m (m-1)(m-2)n_m'.$$
This equality explains the contribution of the points in $A$. Now let $b \in B$ such that $n(b)=m$. The number of such points is precisely $n_m''$. It follows that
$$\sum_{b\in B}(n(b)-1)=\sum_m(m-1)n_m''.$$
Let $1 \leq i<j\leq r$ and consider the intersection $C_i \cap C_j$. It contains exactly $N_iN_j$ points, since $C_i$ and $C_j$ intersects transversely. The sum $S=\sum_{1 \leq i<j\leq r}N_iN_j$
represents the number of all such intersection points. Note that a point  $b\in B$ is counted in this sum exactly ${n(b) \choose 2}$ times. This yields the following formula
$$2S=\sum_mm(m-1)n_m''.$$
These formulas give the following result.

\begin{thm} \label{mhnU2}
With the above assumptions and notation, one has
$$\dim Gr^2_FH^2(U,\C)=\frac{(N-1)(N-2)}{2} -\sum_m{m-1 \choose 2}n_m,$$
with $n_m=n_m'+n_m''$ the number of ordinary $m$-tuple points of $C$.

\end{thm}

The following consequence of Theorem \ref{mhnU} and Theorem \ref{mhnU2} applies in particular to any projective line arrangement.

\begin{cor}\label{linearr}
Assume that   $C=\bigcup_{j=1,r}C_j$ is the decomposition of $C$ as a union of irreducible components $C_j$, with any curve $C_j$ having only ordinary multiple points as singularities and being rational, i.e. $g_j=0$.
If the intersection of any two distinct components $C_i$ and $C_j$ is transverse, i.e. the points in $C_i \cap C_j$ are nodes of the curve $C_i \cup C_j$, then one has 
$$\dim H^2(U,\C)=\frac{(N-1)(N-2)}{2} -\sum_m{m-1 \choose 2}n_m,$$
with $n_m$ the number of ordinary $m$-tuple points of $C$.

\end{cor}

\section{ Curves with ordinary singularities of multiplicity $\leq 4$}

Let $C\subset \mathbb{P}^2$ be a curve of degree $N$ having only ordinary singular points of multiplicity at most $4$. Set $U=\mathbb{P}^2\setminus C$, and let $C=\cup_{j=1}^rC_j$ be the decomposition of $C$ in irreducible components. Then,
\begin{eqnarray*}
P(C)&=&\sum_{j=1}^r P(C_j)-\sum_{0\leq i <j\leq r} P(C_i\cap C_j) + \sum_{0\leq i <j<k\leq r} P(C_i\cap C_j\cap C_k)\\ 
&-& \sum_{0\leq i <j<k<l\leq r} P(C_i\cap C_j\cap C_k \cap C_l).
\end{eqnarray*}
Let $a_m^j$ denote the number of singular points of multiplicity $m$ that belong to the  component $C_j$ (note that a point can be singular on two components, being a node on each of them).\\
 Denote by $b_3^k$ (respectively $b_4^k$) the number of triple points (respectively points of multiplicity 4) of $C$ that are intersection of exactly $k$ components, for $k=2,3$ (respectively $k=3,4$). Let $b_4^2$ (respectively $\tilde{b_4^2}$) be the number of singular points $p$ of multiplicity $4$ in $C$ representing the intersection of exactly 2 components, such that one of which has a triple point at $p$ (respectively each one has a node at $p$). Then one has
$$\sum_{0\leq i <j\leq r} P(C_i\cap C_j)=\sum_{0\leq i <j\leq r} N_iN_j - b^2_3-3\tilde{b_4^2}-2b_4^2-2b_4^3.$$
Indeed, a point of type $b_3^2$ (resp. $b_4^2$, resp. $\tilde{b_4^2}$) occurs only in one intersection $C_i \cap C_j$, and has the multiplicy 2 (resp.3, resp. 4) in this intersection. A point of type $b_4^3$ occurs in 3 intersections $C_i \cap C_j$ with multiplicitities $1,2,2$, and this accounts for the correction term $-2b_4^3$. Then one has
$$\sum_{0\leq i <j<k\leq r} P(C_i\cap C_j\cap C_k)= b^3_3+b_4^3+{4\choose 3}b_4^4,$$
and
$$\sum_{0\leq i <j<k<l\leq r} P(C_i\cap C_j\cap C_k \cap C_l)=b^4_4.$$

Hence, by Proposition \ref{propirred}, we get the following.
\begin{eqnarray*}
P(C)&=&ruv-(\sum_{j=1}^r g_j) u-(\sum_{j=1}^r g_j) v-\sum_{j=1}^r (a_2^j+2a_3^j+3a^j_4)-\sum N_iN_j\\
&+&b_3^2+3\tilde{b^2_4}+2b^2_4+3b_4^3+b^3_3+3b^4_4.
\end{eqnarray*}
Therefore, as above, we obtain
\begin{eqnarray*}
P(U)&=&u^2v^2-(r-1)uv+1-r+(\sum_{j=1}^r g_j) u+(\sum_{j=1}^r g_j) v+\sum_{j=1}^r (a_2^j+3a_3^j+6a^j_4)\\
&-&\sum_{j=1}^r(a_3^j+3a_4^j)+\sum N_iN_j-b_3^2-3\tilde{b^2_4}-2b^2_4-3b_4^3-b^3_3-3b^4_4.
\end{eqnarray*}
Finally we get \begin{eqnarray*}
\dim Gr^2_F H^2(U)&=&\sum_{j=1}^r (g_j+a_2^j+3a_3^j+6a^j_4-1)+\sum N_iN_j+1-(\sum_{j=1}^r a_3^j+b_3^2+b^3_3) \\
&-&3(\sum_{j=1}^r a_4^j+\tilde{b^2_4}+b^2_4+b_4^3+b^4_4)+{b^2_4}\\
&=&\frac{(N-1)(N-2)}{2}-n_3-3n_4+{b_4^2},
\end{eqnarray*}
with $n_m$ the number of ordinary $m$-tuple points of $C$.

\begin{thm} \label{mhnUformult4}
Let $C\subset \mathbb{P}^2$ be a curve of degree $N$ having only ordinary singular points of multiplicity at most $4$. If $U=\mathbb{P}^2\setminus C$,
then one has
$$\dim Gr^2_FH^2(U,\C)=\frac{(N-1)(N-2)}{2} -\sum_{m=3,4}{m-1 \choose 2}n_m +b_4^2,$$
with $n_m$ the number of ordinary $m$-tuple points of $C$ and $b_4^2$ the number of singular points $p$ of $C$ which are smooth on one component $C_i$ of $C$ and have multiplicity $3$ on the other component $C_j$ of $C$ passing through $p$.

\end{thm}

\section{Pole order filtration versus Hodge filtration for plane curve complements}

For any hypersurface $V$ in a projective space $\PP^n$, the cohomology groups $H^*(U,\C)$ of the complement $U =\PP^n \setminus V$ have a pole order filtration $P^k$, see for instance \cite{DSt2}, and it is known by the work of P. Deligne, A. Dimca \cite{DD} and M. Saito \cite{MS}  that one has
$$F^kH^m(U,\C) \subset P^kH^m(U,\C)$$
for any $k$ and any $m$. For $m=0$ and $m=1$, the above inclusions are in fact equalities
(the case $m=0$ is obvious and the case $m=1$ follows from the equality $F^1H^1(U,\C)=H^1(U,\C)$). For $m=2$, we have again $F^kH^2(U,\C) = P^kH^2(U,\C)$ for $k=0,1$ for obvious reasons, but one may get strict inclusions 
$$F^2H^2(U,\C) \ne P^2H^2(U,\C)$$
already in the case when $V=C$ is a plane curve, see \cite{DS}, Remark 2.5 or \cite{Dbk}.
However, to give such examples of plane curves was until now rather complicated. We give below a numerical condition which tells us exactly when the above strict inclusion holds. 

We need first to recall some basic definitions.
Let $S=\oplus_rS_r=\C[x,y,z]$ be the graded ring of polynomials with complex coefficients, where $S_r$ is the vector space of homogeneous polynomials of $S$ of degree $r$. For a homogeneous polynomial $f$ of degree $N$, define the Jacobian ideal of $f$ to be the ideal $J_f$ generated in $S$ by the partial derivatives $f_x,f_y,f_z$ of $f$ with respect to $x$, $y$ and $z$. The graded \textit{Milnor algebra} of $f$ is given by $$M(f)=\oplus_rM(f)_r=S/J_f.$$
Note that the dimensions $\dim M(f)_r$ can be easily computed in a given situation using some computer software e.g. Singular.
Now we can state the main result of this section.

\begin{thm} \label{F=P}
Let $C:f=0$  be a reduced curve of degree $N$ in $\mathbb{P}^2$ having only weighted homogeneous singularities  and let $C_i$ for $i=1,...,r$ be the irreducible components of $C$. If $U=\mathbb{P}^2\setminus C$,
then
$$\dim P^2H^2(U,\C)-\dim F^2H^2(U,\C) =\tau(C) + \sum_{i=1,r}g_i-\dim M(f)_{2N-3},$$
where $\tau(C)$ is the global Tjurina number of $C$ (that is the sum of the Tjurina numbers of all the singularities of $C$) and  $g_i$ is the genus of the normalization of $C_i$ for $i=1,...,r$.
\end{thm}
In particular we get the following result, which yields in particular a new proof for Theorem 1.3 in 
\cite{DSt}.

\begin{cor}\label{rational}
If a reduced plane curve has only nodes as singularities,  then one has
$$\dim M(f)_{2N-3}=\tau(C) + \sum_{i=1,r}g_i.$$
\end{cor}

\proof Indeed, it is known that for a nodal curve one has the equality $F^2H^2(U,\C) = P^2H^2(U,\C)$, see \cite{De} or \cite{MS}.

\endproof
Note that we  have the following obvious consequence of Theorem  \ref{mhnU}.
\begin{cor}\label{stabilization}
For a reduced plane curve  $C$ one has 
$$\dim P^2H^2(U,\C)-\dim F^2H^2(U,\C) \leq  \sum_{i=1,r}g_i.$$

\end{cor}

\proof  Indeed, Theorem  \ref{mhnU} can be restated as 
$$\dim H^2(U,\C)-\dim F^2H^2(U,\C) = \sum_{i=1,r}g_i,$$
in view of the equality $F^1H^2(U,\C)=H^2(U,\C)$, see \cite{Dbk}, proof of Corollary 1.32, page 185.
\endproof

\begin{rk}\label{rkrational}
If a reduced plane curve $C$ has only rational irreducible components, i.e. $g_i=0$ for all $i$, then the above inequality implies $F^2H^2(U,\C) = P^2H^2(U,\C)$. 
This result can be regarded as an improvement of a part of the Remark 2.5 in \cite{DS}, where the result is claimed only for curves with nodes and cusps as singularities. 

\end{rk}
The above discussion implies also the following result, which can be regarded as a generalization of Theorem 4.1 (A) in \cite{A}.

\begin{cor}\label{rational2}
If a reduced plane curve $C:f=0$ has only weighted homogeneous  singularities,  then one has
$$0\leq \dim M(f)_{2N-3}-\tau(C) \leq  \sum_{i=1,r}g_i.$$
In particular, if in addition  the curve $C$ has only rational irreducible components,  then one has
$$\dim M(f)_{2N-3}=\tau(C) .$$
\end{cor}

Now we give the proof of Theorem \ref{F=P}.
Corollary 1.3 in \cite{DSt2} implies that
$$\dim P^2H^2(U,\C)= \dim H^2(U,\C)+\tau(C)-\dim M(f)_{2N-3}.$$ 
On the other hand, Theorem \ref{mhnU} and the fact $\dim F^1H^2(U,\C)=  H^2(U,\C)$ yield
$$\dim F^2H^2(U,\C)= \dim H^2(U,\C)-\sum_{i=1,r}g_i,$$
which clearly completes the proof of Theorem \ref{F=P}.

\begin{ex}\label{free}
In this example we present a free divisor $C:f=0$, whose irreducible components consist of 12 lines and one elliptic curve, and where $F^2H^2(U,\C) \ne P^2H^2(U,\C)$.
Let $f=xyz(x^3+y^3+z^3)[(x^3+y^3+z^3)^3-27x^3y^3z^3].$ If we consider the pencil of cubic curves $(x^3+y^3+z^3, xyz)$, then the curve $C$ contains all the singular fibers of this pencil, and this accounts for the 12 lines given by 
$$xyz[(x^3+y^3+z^3)^3-27x^3y^3z^3]=0,$$
and the elliptic curve (hence of genus 1) given by $x^3+y^3+z^3=0$.
Then $C$ is a free divisor, see \cite{JV} or by a direct computation using Singular, which shows that $I=J_f$, where $I$ is the saturation of the Jacobian ideal $J_f$, see Remark 4.7 in \cite{DSe}. The direct computation by Singular also yields $\tau(C)=156$ and $\dim M(f)_{2N-3}=\dim M(f)_{27}=156$. Moreover, applying Corollary 1.5 in \cite{DSt3}, we see via a Singular computation that all singularities of the curve $C$ are weighted homogeneous. Alternatively, there are 12 nodes, 3 in each of the 4 singular fibers of the pencils (which are triangles), and the 9 base points of the pencil, each an ordinary point of multiplicity 5. Each of the 12 lines contains exactly 3 of these base points, and they are exactly the intersection of the elliptic curve with the line. This description implies that there are no other singularities, in accord with
$$12+9 \times 16=156= \tau(C).$$
It follows from
 Theorem \ref{F=P} that
$\dim P^2H^2(U,\C)-\dim F^2H^2(U,\C) =1.$
Hence the presence of a single irrational component of $C$ leads to $F^2H^2(U,\C) \ne P^2H^2(U,\C)$.
\end{ex}

\bigskip

\small \textbf{Acknowledegment:} I gratefully acknowledge the support of the Lebanese National Council for Scientific Research, without which the present study could not have been completed.

\end{document}